%% file: genspatimparpod.tex
\documentclass[a4paper,10pt]{article}
\usepackage[top=90pt,bottom=80pt,left=100pt,right=100pt]{geometry}
\usepackage[T1]{fontenc}
\usepackage[utf8]{inputenc}
\usepackage{lmodern}
\usepackage{epsfig}
\usepackage{amssymb,amsmath,amsthm,mathrsfs,mathtools,amscd}
\usepackage{float}
\usepackage[figtopcap]{subfigure}
\usepackage{color}
\usepackage{authblk}
\usepackage{tikz}
\usepackage{pgfplots}
\usepackage{framed}
\usepgfplotslibrary{groupplots}
\pgfplotsset{compat=newest}
\usepackage{subfigure}

\usepackage[bookmarks=true]{hyperref}

\input{stylefyle} %uses mathcomabb.tex to define commands for symbols and notations 
\input{mathComAbb} %uses mathcomabb.tex to define commands for symbols and notations 

\title{\bf Space-time Galerkin POD with application in optimal control of semi-linear parabolic partial differential equations}
\author[1]{\href{mailto:m.m.baumann@tudelft.nl}{Manuel Baumann}}
\author[2]{\href{mailto:benner@mpi-magdeburg.mpg.de}{Peter Benner}}
\author[2]{\href{mailto:heiland@mpi-magdeburg.mpg.de}{Jan Heiland\footnote{Corresponding author's e-mail address: \texttt{heiland@mpi-magdeburg.mpg.de}}}}

\affil[1]{Delft Institute of Applied Mathematics}
\affil[2]{Max Planck Institute for Dynamics of Complex Technical Systems Magdeburg}
\begin{document}
\maketitle
% \tableofcontents

\begin{abstract}
In the context of Galerkin discretizations of a partial differential equation (PDE), the modes of the classical method of Proper Orthogonal Decomposition (POD) can be interpreted as the ansatz and trial functions of a low-dimensional Galerkin scheme. If one also considers a Galerkin method for the time integration, one can similarly define a POD reduction of the temporal component. This has been described earlier but not expanded upon -- probably because the reduced time discretization globalizes time which is computationally inefficient. However, in finite-time optimal control systems, time \textit{is} a global variable and there is no disadvantage from using a POD reduced Galerkin scheme in time.
In this paper, we provide a newly developed generalized theory for space-time Galerkin POD, prove its optimality in the relevant function spaces, show its application for the optimal control of nonlinear PDEs, and, by means of a numerical example with Burgers' equation, discuss the competitiveness by comparing to standard approaches.
\end{abstract}

\input{mathEnv}

\input{chap-introduction}
\input{chap-spacetimegenpod}

\input{chap-galerkinpodspacetime}
\input{chap-implementationissues}

\input{chap-optcontrol}

\input{chap-numerics}

\input{chap-conclusion}

\bibliographystyle{abbrv}
% \bibliography{genspactimpod}
% \bibliography{csc,bib_jh,swbib_jh}

\input{bib.bbl}
% \bibliography{genspactimpod,bib_mb,bib_jh,mor}
% \bibliography{bib_mb,bib_jh,mor,swbib_jh,genspactimpod}

\end{document}

%% file: stylefyle.tex
\def\trp{{\mathsf{T}}}  % transpose should be different from t_e

\def\ouss{Y}
\def\oust{S}

\def\Ouss{\mathcal{\ouss}}
\def\Oust{\mathcal{\oust}}
\def\hOuss{\hat{\mathcal{\ouss}}}
\def\hOust{\hat{\mathcal{\oust}}}
\def\Ousts{\ensuremath{{\Oust \cdot \Ouss}}}
\def\hOusts{\ensuremath{{\hOust \cdot \hOuss}}}

\def\Somts{\ensuremath{\mathcal{X}}}

\def\hAdjss{\hat \Lambda}
\def\Adjst{{\mathcal R}}
\def\hAdjst{\hat \Adjst}

\def\Hoi{H_0^1(\Omega)}
\def\Hiot{H^1(0, T)}

\newcommand{\stdp}[2]{((#1,#2))_{\Oust\cdot\Ouss}}

\def\xtsvec{\ensuremath{\mathbf x}}
\def\hxtsvec{\ensuremath{\hat{\mathbf x}}}
\newcommand{\xtsvecxx}[2]{\xtsvec_{#1\cdot#2}}
\newcommand{\hxtsvecxx}[2]{\hxtsvec_{#1\cdot#2}}
\def\xtsvecij{\ensuremath{\mathbf x_{i\cdot j}}}
\def\hUpsilon{{\hat \Upsilon}}
\def\hPsi{{\hat \Psi}}

\def\My{\mathbf{M}_\Ouss}
\def\Ms{\mathbf{M}_\Oust}
\def\Ly{\mathbf{L}_\Ouss}
\def\Ls{\mathbf{L}_\Oust}
\def\Lyt{\Ly^\trp}
\def\Lymt{\Ly^{-\trp}}
\def\Lst{\Ls^\trp}

\def\brgv{z}

\def\hnu{\hat \nu}
\def\hpsi{\hat \psi}
\def\hdms{d \hms}
\def\hms{M _{\hat{\mathcal S}}}
\def\hmy{M _{\hat{\mathcal Y}}}
\def\hay{K _{\hat{\mathcal Y}}}
\def\hbv{\hat{\mathbf v}}
\def\hnsy{H_{\hat{\mathcal S}\hat{\mathcal Y}}}
\def\hfsy{f_{\hat{\mathcal S}\hat{\mathcal Y}}}
\def\hsbas{\hat \Psi}
\def\hybas{\hat \Upsilon}

\def\ybas{\Upsilon}

\def\hats{\hat s}
\def\hatr{\hat r}
\def\hatq{\hat q}
\def\hatp{\hat p}

\def\hphi{\hat \phi}
\def\hmu{\hat \mu}
\def\hml{M _{\hAdjss}}
\def\hmr{M _{\hAdjst}}
\def\hdmr{dM _{\hAdjst}}
\def\hal{K _{\hAdjss}}

\def\hadjvd{\hat{\boldsymbol \lambda}}
\def\hnrl{D_x N_{\hAdjss\hAdjst}^\trp}
\def\hmyl{M _{\hat{\mathcal Y}\hAdjss}}
\def\hmly{M _{\hAdjss\hat{\mathcal Y}}}
\def\hmsr{M _{\hOust \hAdjst}}
\def\hmrs{M _{\hAdjst\hOust}}

\def\nsnap{\ensuremath{s}}

\providecommand{\bbmat}{\begin{bmatrix}}
\providecommand{\ebmat}{\end{bmatrix}}

\def\tsint{\int_0^T \int_\Omega}
\def\tint{\int_0^T}
\def\sint{\int_\Omega}
\providecommand{\partielt}[1]{\frac{\partial #1}{\partial t}}

\def\ctime{\texttt{walltime}}
\def\vterm{\ensuremath{{\frac 12 \norm{\hat x - x_0}_{L^2}^2}}}
\def\cfune{\ensuremath{\mathcal J(\hat x, \hat u)}}

%% file: mathComAbb.tex
\DeclareMathOperator{\spann}{span}

\DeclareMathOperator{\mattovec}{vec}

\def\trp{{\mathsf T}}

\providecommand{\norm}[1]{\lVert#1 \rVert}

\providecommand{\inva}[1]{\text{~\textup{d}} #1}
\providecommand{\andi}[0]{\quad \text{and} \quad}

\providecommand{\partielt}[1]{\frac{\partial #1}{\partial t}}

\providecommand{\bbmat}{\begin{bmatrix}}
\providecommand{\ebmat}{\end{bmatrix}}

%% file: mathEnv.tex
\newtheorem{thm}{Theorem}[section]
\newtheorem{cor}[thm]{Corollary}
\newtheorem{lem}[thm]{Lemma}
\newtheorem{ass}[thm]{Assumption}

\theoremstyle{remark}
\newtheorem{exa}[thm]{Example}

\theoremstyle{remark}
\newtheorem{rem}[thm]{Remark}

\theoremstyle{definition}
\newtheorem{defin}[thm]{Definition}

\theoremstyle{plain}
\newtheorem{prob}[thm]{Problem}

\newtheoremstyle{named}{}{}{\itshape}{}{\bfseries}{.}{.5em}{#1 \thmnote{#3}}
\theoremstyle{named}
\newtheorem*{namedproblem}{Problem}

\theoremstyle{plain}
\newtheorem{prop}[thm]{Proposition}

\theoremstyle{plain}
\newtheorem{alg}[thm]{Algorithm}

\floatstyle{ruled}
\newfloat{algorithm}{h}{loa}[section]
\floatname{algorithm}{Algorithm}

%% file: chap-introduction.tex
\section{Introduction}

The method of \emph{Proper Orthogonal Decomposition} (POD) is a standard model reduction tool. For a generic dynamical system
\begin{equation}\label{eq:gendynsys}
  \dot v = f(t, v),
\end{equation}
on the time interval $(0, T]$ with a solution $v$ with $v(t) \in \mathbb R^{N}$ and using samples $v(t_j)$, POD provides a set of $\hat n$ so-called modes $\hat v_1$, ..., $\hat v_{\hat n} \in \mathbb R^{N}$ which optimally parametrize the solution trajectory. As a result, the system \eqref{eq:gendynsys} can be projected down to a system of reduced spatial dimension $\hat n$ that reflects the dynamical behavior of \eqref{eq:gendynsys} well.
If the considered system stems from a \emph{Finite Element} (FEM) discretization of a PDE, then the modes $\hat v_i$, $i=1,\dotsc,\hat n$ can be interpreted as ansatz functions in the finite element space $\Ouss$ and the projected system as a particular Galerkin projection of the underlying PDE.

In this paper we provide a theoretical framework and show cases for a space-time Galerkin POD method. The underlying ideas for this generalization of POD have been developed and tested in our earlier works \cite{BauBH15,BauHS15}.

The first innovation of the proposed generalized POD approach bases on the observation that instead of the discrete time samples $v(t_j)$ one may use the projection of $v$ onto the finite dimensional subspace $\Oust\cdot \Ouss$, where $\Oust$ is a, say, $k$-dimensional subspace of $L^{2}(0,T)$. The second innovation is that the projection onto $\Ousts$ can be interpreted as Galerkin discretization in time which can be reduced analoguously to the POD reduction of the space dimension. The resulting scheme is a POD reduced space-time Galerkin discretization.

This basic idea of a space-time POD has already been taken up in \cite{VolW06} but not progressed since then. We think that this is due to the fact that temporal POD destructs the causality in time which makes it very inefficient for numerical simulations. In fact, the POD reduced time ansatz functions are global such that the space-time Galerkin system has to be solved as a whole rather than in sequences of time slobs as in standard time-stepping or discontinuous Galerkin schemes \cite{MDGC16,Zw14}. Thus, the reduced space-time scheme cannot compete with, e.g., a spatial POD combined with a standard Runge-Kutta solver. However, in finite-time optimal control problems, the time \emph{is} a global variable and, as we will show by numerical examples, the space-time Galerkin discretization becomes very competitive.

The need and the potential of also reducing the time dimension of a reduced order model have been discussed in \cite{CarRB15}. There -- similar to our observation that an SVD of a matrix of measurements also reveals compressed time information -- it is proposed to use the right singular vectors of a classical snapshot matrix for forecasting.

We want to point out that the method of \emph{Proper Generalized Decomposition} (PGD) is related to the proposed space-time Galerkin POD only in so far as for PGD also space-time (and parameter) tensor bases are used for the modelling; see, e.g., \cite{ChiALK11}. However, the PGD approach seeks to successively build up the bases by collocation, \emph{greedy algorithms}, and fixed-point iteration, whereas our approach reduces a given basis on the base of measurements. For the same reasons, the connection of the presented approach to other tensor-based low-dimensional approximation schemes \cite{KhoS11, SchS09} as well as to \emph{Reduced Basis} approaches \cite{YanPU14} is only marginal.

This paper is organized as follows. At first, we introduce the mathematical framework and rigorously prove the optimality of the reduced space and time bases. Then we illustrate how the reduced bases can be used for low-dimensional space-time Galerkin approximations. In particular, we address how to treat quadratic nonlinearities, how to incorporate initial and terminal values, and how to set up the bases for a general PDE by means of standard approximation schemes. Finally, we illustrate the performance of the space-time Galerkin POD approach for the optimal control of Burgers' equation and compare it to well-established gradient-based methods combined with standard POD.

%% file: chap-spacetimegenpod.tex
\section{Space-Time Galerkin POD}
In this section, we provide the analytical framework for space-time POD. We introduce the considered function spaces and directly prove the optimality of the POD projection in the respective space-time $L^2$ norm. For a time interval $(0,T)$ and a spatial domain $\Omega$, consider the space-time function space $L^2( 0,T;L^2( \Omega))$. Let 
\begin{equation*}
	\Oust=\spann\{\psi_1, \dotsc, \psi_s\}\subset L^2( 0,T)\quad\text{and}\quad
	\Ouss = \spann\{\nu_1, \dotsc, \nu_q\}\subset L^2( \Omega)
\end{equation*}
be finite dimensional subspaces of dimension $s$ and $q$, respectively, and let 
\begin{equation} \label{eq:spaceofmeas}
	\Somts = \mathcal S\cdot \mathcal Y \subset L^2( 0,T;L^2( \Omega)).
\end{equation}

The space-time $L^2$-orthogonal projection $\bar x := \Pi_\Ousts x$ of a function $x\in L^2(0,T; L^2(\Omega))$ onto \Somts~is given as

\begin{equation}
	\bar x (\xi,\tau) = \sum_{j=1}^s \sum_{i=1}^q \xtsvecij \nu_i(\xi) \psi_j(\tau),
\end{equation}
where the coefficient $\xtsvecij$ are the entries of the matrix
\begin{equation} \label{eq:somprojmatrix}
	\mathbf{X} = \bigl[ \xtsvecij \bigr]_{i=1,\dotsc,q}^{j=1,\dotsc,s}:=\My^{-1}
	\begin{bmatrix} 
		\stdp{x}{\nu_1\psi_1} & \hdots & \stdp{x}{\nu_1\psi_s} \\ \vdots & \ddots & \vdots \\ \stdp{x}{\nu_q\psi_1} & \hdots & \stdp{x}{\nu_q\psi_s} 
	\end{bmatrix}
	\Ms^{-1},
	%\in \mathbb{R}^{q \times s},
\end{equation}
where 
$$\stdp{x}{\nu_i\psi_j}:= ( (x,\nu_i)_\Ouss, \psi_j)_\Oust:=\int_0^T\bigr( \int_\Omega x(\xi, \tau)\nu_i(\xi)\inva \xi\bigr) \psi_j(\tau) \inva \tau.$$

Here, $\My^{-1}$ and $\Ms^{-1}$ are the inverses of the mass matrices with respect to space and time,
\begin{equation} \label{eq:massmatrices}
\My := \bigl[ (\nu_i, \nu_j)_\Ouss \bigr]_{i=1,\dotsc,q}^{j=1,\dotsc,q} \quad \text{and} \quad \Ms := \bigl[ (\psi_i, \psi_j)_\Oust \bigr]_{i=1,\dotsc,s}^{j=1,\dotsc,s}.
\end{equation}

\begin{rem}\label{rem:measuremeasure}
	We will refer to $\mathcal X = \Ousts$ as the measurement space, to the basis functions of $\Ouss$ and $\Oust$ as measurement functions, and to $\mathbf X$ as the measurement matrix. This means that a function in $L^2( 0,T;L^2( \Omega))$ can be measured in $\mathcal X$, e.g. via its projection on $\mathcal X$, and, the other way around, an element $\mathbf X$ of $\mathcal X$ can be seen as a measurement of some functions in $L^2( 0,T;L^2( \Omega))$.
\end{rem}

We introduce some representations of the inner product and the norm of functions in $\Ousts$.
\begin{lem}[Space-time discrete $L^2$-product]\label{lem:st-disc-prod-norm}
Let
\begin{align*}
x^1 = \sum_{j=1}^s \sum_{i=1}^q \xtsvecij^1 \nu_i \psi_j \in \Ousts, \quad
x^2 = \sum_{j=1}^s \sum_{i=1}^q \xtsvecij^2 \nu_i \psi_j \in \Ousts,
\end{align*}
then, with 
\begin{align*}
	\xtsvec^\ell = [\xtsvecxx 11^\ell,\dotsc,\xtsvecxx q1^\ell,
									\xtsvecxx 12^\ell,\dotsc,\xtsvecxx q2^\ell,
									\dotsc,
									\xtsvecxx 1s ^\ell,\dotsc,\xtsvecxx qs^\ell]^\trp 
				=: \mattovec(\mathbf{X^\ell}), \quad \ell = 1, 2,
\end{align*}
the inner product in \Ousts~is given as
\begin{equation}
	\stdp{x^1}{x^2}= \int_0^T\int_\Omega x^1x^2\inva \xi \inva \tau = (\xtsvec^1)^\trp\left(\Ms \otimes \My \right) \xtsvec^2
	\label{eq:tsvecinnerproduct}
\end{equation}
and the induced norm as
\begin{equation}\label{eq:tsvecnorm}
	\norm{x^\ell}^2_{\Ousts} := \stdp{x^\ell}{x^\ell}= \norm{\xtsvec^\ell}_{\Ms \otimes \My}^2 = \norm{\My^{1/2} \mathbf X^\ell \Ms^{1/2}}_F^2, \quad \ell=1,2,
\end{equation}
where $\norm{\cdot}_{\Ms \otimes \My}$ denotes the Euclidean vector norm weighted by $\Ms \otimes \My$.
\end{lem}
\begin{proof} Straight-forward calculations. \end{proof}
\begin{rem}\label{rem:cholesky}
In practical applications, one uses a \emph{Cholesky}-factorization of the mass matrices \eqref{eq:massmatrices} rather than the square-root.
\end{rem}
\begin{cor}\label{cor:norminchols}
	Let $\Ms=\Ls\Lst$ and $\My=\Ly\Lyt$ be given in factored form. Then, for a given $x \in \Ousts$ with its coefficient matrix $\mathbf X$ and vector $\xtsvec = \mattovec{(\mathbf X)}$ it holds that
\begin{equation}
	\norm{x}^2_{\Ousts} = \norm{\xtsvec}_{\Ms \otimes \My}^2 = \norm{\Lyt \mathbf X\Ls}_F^2.
\end{equation}
\end{cor}
\begin{proof}
	\begin{align*}
		\norm{\xtsvec} _{\Ms \otimes \My}^2 
		& = \xtsvec^\trp (\Ms \otimes \My ) \xtsvec 
		= \xtsvec^\trp ( \Ls \otimes \Ly ) \cdot (\Lst \otimes \Lyt) \xtsvec \\
		& = \norm{ (\Lst \otimes \Lyt ) \xtsvec }_2^2 
		= \norm{ \mattovec(\Lyt \mathbf X \Ls) }_2^2 
		= \norm{ \Lyt \mathbf X \Ls }_F^2,
	\end{align*}
  as it follows from basic properties and relations between the Kronecker product, the vectorization operator, and the Frobenius norm.
\end{proof}

From now on, we will always consider the factorized form. In theory, one can always replace the factors by the square roots of the respective mass matrices.

Next, we will consider a given function $x\in\Ousts$ and determine low-dimensional subspaces of $\Ouss$ and $\Oust$ that can provide low-dimensional approximations to $x$ in a norm-optimal way.

\begin{lem}[Optimal low-rank bases in space]\label{lem:optspacebas}
	Given $x\in\Ousts$ and the associated matrix of coefficients $\mathbf X$. The best-approximating $\hat q$-dimensional subspace $\hat \Ouss$ in the sense that 
	$\norm{x-\Pi_{\Oust\cdot \hat \Ouss}x}_\Ousts$ 
	is minimal over all subspaces of $\Ouss$ of dimension $\hat q$ is given as $\spann\{\hat \nu_i\}_{i=1,\dotsc,\hat q}$, where 
	\begin{equation}
\begin{bmatrix} \hat \nu_1 \\ \hat \nu_2 \\ \vdots \\ \hat \nu_{\hat q} \end{bmatrix} 
	= V_{\hat q}^\trp \Ly^{-1} 
\begin{bmatrix}  \nu_1 \\  \nu_2 \\ \vdots \\  \nu_{q} \end{bmatrix} ,
		\label{eq:opttildebasy}
	\end{equation}
	where $V_{\hat q}$ is the matrix of the $\hat q$ leading left singular vectors of the matrix
	\begin{equation*}
		\Lyt \mathbf X \Ls.
	\end{equation*}
\end{lem}

\begin{proof}
	For the time dimension at fixed index $j$, we consider
	$$y:=\sum_{i=1}^q \xtsvecij \nu_i =
	\begin{bmatrix} 
	  \xtsvecxx 1j &\hdots &\xtsvecxx qj 
	\end{bmatrix} 
	\begin{bmatrix} 
	  \nu_1 \\ \vdots \\ \nu_q 
	\end{bmatrix}
		\in \Ouss .$$

 Next, we determine the orthogonal projection of $y$ onto $\hOuss$. Therefore, we write $y$ as a function in $\hOuss$ and a reminder $\hat R$ in the orthogonal complement:
\begin{equation*}
y = 
\begin{bmatrix} \xtsvecxx 1j &\hdots &\xtsvecxx qj \end{bmatrix} 
	\begin{bmatrix} \nu_1 \\ \vdots \\ \nu_q \end{bmatrix}
		= 
		\begin{bmatrix} \beta_1 &\hdots &\beta_{\hat q} \end{bmatrix} 
		\begin{bmatrix}\hat  \nu_1 \\ \vdots \\ \hat \nu_{\hat q} \end{bmatrix}
			+ \hat R.
\end{equation*}
We determine the coefficients $\beta_k$, $k=1,\dotsc,\hat q$ by testing against the basis functions of $\hOuss$. By mutual orthogonality of $\hat \nu_i$, $i=1,\dotsc,\hat q$ and their orthogonality against $\hat R$, it follows that
\begin{align*}
  \beta_k = (\sum_{i=1}^{\hat q}\beta_i\hat \nu_i,\hat  \nu_k)_\Ouss = (\hat R + \sum_{i=1}^{\hat q}\beta_i\hat \nu_i,\hat  \nu_k)_\Ouss &= (\sum_{i=1}^{q}\xtsvecxx ij\nu_i,\hat  \nu_k)_\Ouss \\
	&\overset{(*)}{=}
	\begin{bmatrix} 
	  \xtsvecxx ij &\hdots &\xtsvecxx qj 
	\end{bmatrix} 
	% \int_\Omega
	% \begin{bmatrix} 
	%   \nu_1 \\ \vdots \\ \nu_q 
	% \end{bmatrix}
	% \begin{bmatrix} 
	%   \nu_1 & \hdots & \nu_q 
	% \end{bmatrix}
	% \inva \xi
	\My
	\Ly^{-\trp}V_{\hat q,k} \\
  &\overset{\phantom{(*)}}{=} 
	\begin{bmatrix} 
	  \xtsvecxx ij &\hdots &\xtsvecxx qj 
	\end{bmatrix} 
	\Ly V_{\hat q,k},
\end{align*}
where in $\overset{(*)}{=}$ we have used that $ \hat  \nu_k = \begin{bmatrix} \nu_1 & \hdots&  \nu_q \end{bmatrix} \Ly^{-\trp}V_{\hat q, k}$ and where $ V_{\hat q, k}$ is the $k$-th column of $V_{\hat q}$ in \eqref{eq:opttildebasy}.  Thus, we find that the coefficients of the orthogonal projection of $y$ onto $\hat \Ouss$ in the bases of $\hOuss$ and $\Ouss$ are given through 
	\begin{align*}
		\hat y=\sum_{i=1}^{\hat q}\beta_i\hat\nu_i =
		\begin{bmatrix} \beta_1 &\hdots &\beta_{\hat q} \end{bmatrix} 
		\begin{bmatrix}\hat  \nu_1 \\ \vdots \\ \hat \nu_{\hat q} \end{bmatrix}
			&= \begin{bmatrix} \xtsvecxx 1j &\hdots &\xtsvecxx qj \end{bmatrix}
				\Ly V_{\hat q} \begin{bmatrix} \hat \nu_1 \\ \vdots \\ \hat  \nu_q \end{bmatrix} \\
			&= \begin{bmatrix} \xtsvecxx 1j &\hdots &\xtsvecxx qj \end{bmatrix}
				\Ly V_{\hat q}V_{\hat q}^\trp \Ly^{-1}
			\begin{bmatrix} \nu_1 \\ \vdots \\ \nu_q \end{bmatrix} \\
				&=: \begin{bmatrix}{\hxtsvecxx 1j} & \hdots & {\hxtsvecxx qj} \end{bmatrix}
				\begin{bmatrix} \nu_1 \\ \vdots \\ \nu_q \end{bmatrix}.
		\end{align*}
		Noting that 
		$
	\begin{bmatrix} 
	  \xtsvecxx 1j &\hdots &\xtsvecxx qj 
	\end{bmatrix} ^{\trp}
	$ makes up the $j$-th column of the matrix $\mathbf X$ associated with $x$, we conclude that the matrix $\hat {\mathbf X}$ of coefficients associated with $\Pi_{\Oust\cdot \hat \Ouss}x$ is given as 
		$$\hat {\mathbf X} = \Lymt V_{\hat q}V_{\hat q}^\trp \Lyt \mathbf X$$ 
		and, by Corollary \ref{cor:norminchols}, we have that
		\begin{align*}
			\norm{x - \Pi_{\Oust\cdot \hat \Ouss}x}_\Ousts
			& =\norm{\Lyt{\mathbf X}\Ls-\Lyt \hat{ \mathbf X} \Ls}_F 
			= \norm{\Lyt [{\mathbf X}-\hat{ \mathbf X}] \Ls}_F \\
			& = \norm{\Lyt\mathbf X \Ls- V_{\hat q}V_{\hat q}^\trp \Lyt \mathbf X \Ls}_F
		\end{align*}
		which is minimized over all $V_{\hat q}\in \mathbb R^{q,\hat q}$ matrices by taking $V_{\hat q}$ as the matrix of the $\hat q$ leading left singular vectors of $\Lyt \mathbf X \Ls$.
\end{proof}

The same arguments apply to the transpose of $\mathbf X$:
\begin{lem}[Optimal low-rank bases in time]\label{lem:opttimbas}
	Given $x\in\Ousts$ and the associated matrix of coefficients $\mathbf X$. The best-approximating $\hat s$-dimensional subspace $\hat \Oust$ in the sense that 
	$\norm{x - \Pi_{\hat \Oust\cdot \Ouss}x}_\Ousts$ 
	is minimal over all subspaces of $\Oust$ of dimension $\hat s$ is given as $\spann\{\hat \psi_j\}_{j=1,\dotsc,\hat s}$, where 
	\begin{equation}
\begin{bmatrix} \hat \psi_1 \\ \hat \psi_2 \\ \vdots \\ \hat \psi_{\hat s} \end{bmatrix} 
	= U_{\hat s}^\trp \Ls ^{-1}
\begin{bmatrix}  \psi_1 \\  \psi_2 \\ \vdots \\  \psi_{s} \end{bmatrix} ,
		\label{eq:opttildebass}
	\end{equation}
	where $U_{\hat s}$ is the matrix of the $\hat s$ leading right singular vectors of 
	\begin{equation*}
		\Lyt\mathbf X \Ls.
	\end{equation*}
\end{lem}

\begin{rem}\label{rem:podisdegengenpod}
  The approximation results Lemma \ref{lem:optspacebas} and Lemma \ref{lem:opttimbas} hold in the space-time $L^2$ norm, which is the appropriate norm for the considered functions and which is not part of the standard POD approach. However, the need for the right norms have been accounted for through the use of weighted inner products or weighted sums. If one lets $\Oust$ degenerate to a set of Dirac deltas, then Lemma \ref{lem:optspacebas} reduces to the optimality result \cite[Thm. 1.8]{Vol11} for the standard POD approximation in the case that the inner product is weighted with the FEM mass matrix. If one chooses $\Oust$ such that the induced time Galerkin scheme resembles a time discretization by the trapezoidal rule (in fact, for any Runge-Kutta scheme and choice of discretization points there exists a corresponding (discontinuous) Galerkin scheme), then Lemma \ref{lem:optspacebas} reduces to the optimality conditions for the \emph{continuous POD} approach given in \cite[Sec. 1.3]{Vol11}.
\end{rem}

\begin{rem}\label{rem:genpodonlyinspace}
 The idea of generalized measurements also works as a generalization of POD for model order reduction in space. Consider the dynamical system \eqref{eq:gendynsys}, and define $\mathbf{X}_\mathcal{S}:=[(v_i,\psi_j)_\mathcal{S}]_{i=1,...,q}^{j=1,...,s}$, where $v_i$ is the $i$-th component of the vector-valued solution. Then the leading left singular vectors of the matrix $\mathbf{X}_\mathcal{S} \mathbf{L}_\mathcal{S}^{-1}$ are generalized POD modes and a projection of \eqref{eq:gendynsys} onto the space spanned by those modes yields a POD-reduced dynamical system as we have previously described it under the term \emph{gmPOD} in \cite{BauHS15}.
\end{rem}

%% file: chap-galerkinpodspacetime.tex
\section{Space-Time Galerkin Schemes}
In this section, we briefly describe how to formulate a general space-time Galerkin approximation to a generic PDE.
This regression is then followed by the discussion of low-rank space-time Galerkin schemes on the base of POD reductions of standard Galerkin bases.

Let $\{\hat \psi_1, \dotsc, \hat \psi_{\hat s}\}\subset \Hiot$ and $\{\hat \nu_1, \dotsc, \hat \nu_{\hat q}\}\subset \Hoi$ be the POD bases in space and time, respectively. Then, a space-time Galerkin approximation of the generic equation system
\begin{subequations}
	\begin{align}\label{eq:genPDE}
		\dot v - \Delta v + N(v) &= f \quad \hspace{0.12cm} \text{on } (0,T] \times \Omega, \\
		v\bigr|_{\partial \Omega} &= 0 \quad \hspace{0.14cm} \text{on } (0,T],\\
		v\bigr|_{t=0} &= v_0 \quad \text{on } \Omega, \label{eq:genPDE_init}
	\end{align}
\end{subequations}
is given as follows:

The approximate solution $\hat v$ is assumed in the product space $\hat{\mathcal S} \cdot \hat{\mathcal Y}:=\spann\{\hpsi_j\hnu_i\}_{i=1,\dots,\hat q}^{j=1,\dots,\hat s}$. We introduce the formal vectors of the coefficient functions
\begin{equation*}
	\hUpsilon := 
	\begin{bmatrix}
		\hnu_1 \\ \vdots \\ \hnu_{\hat q}
	\end{bmatrix}
	\quad\text{and}\quad
	\hPsi := 
	\begin{bmatrix}
		\hpsi_1 \\ \vdots \\ \hpsi_{\hat s}
	\end{bmatrix}
\end{equation*}
and write $\hat v$ as
\begin{equation}\label{eq:tspodsolansatz}
	% \hat{\mathcal S} \cdot \hat{\mathcal Y} \ni \hat v = 
	\begin{bmatrix}
		\hpsi_1 & \hdots & \hpsi_{\hat q}
	\end{bmatrix}
	\otimes
	\begin{bmatrix}
		\hnu_1 & \hdots & \hnu_{\hat s}
	\end{bmatrix} \hat{\mathbf v}
	=[\hat \Psi^\trp\otimes\hat  \Upsilon^\trp ]\hbv,
\end{equation}
where $\hat{\mathbf v} \in \mathbb R^{\hat s \hat q}$ is the vector of coefficients.
%with respect to the basis of $\hat{\mathcal S}\cdot \hat{\mathcal Y}$. The basis is defined by the outer product of the bases of $\hat{\mathcal S}$ and $\hat{\mathcal Y}$. 
We determine the coefficients by requiring them to satisfy the Galerkin projection of \eqref{eq:genPDE} for every basis function $\hnu_i\hpsi_j$, $i=1,\dots,\hat q$, $j=1,\dots,\hat s$
\begin{equation*}
	\int_0^T \int_\Omega \hnu_i\hpsi_j \dot {\hat v} + \hpsi_j \nabla \hnu_i \nabla \hat v + \hnu_i\hpsi_j N(\hat v)  \inva x \inva t= \int_0^T\int_\Omega \hnu_i\hpsi_j f \inva x \inva t.
\end{equation*}

The latter equations combined give a possibly nonlinear equation system for the vector $\hbv$ of coefficients, which is assembled as follows: For the term with the time derivative we compute
\begin{align*}
	\tsint [\hPsi\otimes \hUpsilon] \partielt{\hat v} \inva x \inva t 
	&= \tsint [\hPsi\otimes \hUpsilon] [\partielt{\hPsi^\trp} \otimes \hUpsilon^\trp] \hbv \inva x \inva t\\
	&= \tsint [\hPsi\partielt{\hPsi^\trp}\otimes\hUpsilon \hUpsilon^\trp] \inva x \inva t \hbv \\
	&= \bigl[\tint \hPsi\partielt{\hPsi^\trp}\inva t \otimes \sint \hUpsilon \hUpsilon^\trp \inva x \bigr]\hbv =: [\hdms \otimes \hmy]\hbv.
\end{align*}
By the same principles, for the term with the spatial derivatives, we obtain
\begin{align*}
	\tsint [\hPsi\otimes \nabla \hUpsilon] \nabla \hat v \inva x \inva t &=
	\tsint [\hPsi\otimes \nabla \hUpsilon] [\hPsi^\trp\otimes \nabla \hUpsilon^\trp]\hbv \inva x \inva t \\
	&= \bigl[\tint \hPsi\hPsi^T \inva t \otimes \sint \nabla \hUpsilon \nabla \hUpsilon^\trp \inva x \bigr ]\hbv := [\hms \otimes \hay ]\hbv.
\end{align*}
Note that in higher spatial dimensions $\nabla \hat v$ as well as $\nabla \hnu_i$ is a vector and, thus, in the preceding derivation, $\nabla \hUpsilon$ has to be interpreted properly.

Summing up, we can write the overall system as

\begin{equation} \label{eq:spatimgalprojsys}
[\hdms \otimes \hmy +\hms \otimes \hay ]\hbv + \hnsy(\hbv) = \hfsy,
\end{equation}
where 
\begin{subequations}\label{eq:spacetimecoeffs}
	\begin{align}
		\hms &:= [\bigl(\hnu_i, \hnu_j \bigr)]_{i,j = 1, \dotsc, \hat s}, \label{eq:timemassmatdirect}\\
		\hdms &:= [\bigl(\hnu_i, \dot \hnu_j \bigr)]_{i,j = 1, \dotsc, \hat s}, \label{eq:timestiffmatdirect}\\
		\hmy &:= [\bigl(\hpsi_l, \hpsi_k \bigr)]_{l,k = 1, \dotsc, \hat q}, \\
		\hay &:= [\bigl(\nabla \hpsi_l, \nabla \hpsi_k \bigr)]_{l,k = 1, \dotsc, \hat q}, \label{eq:spacstiffmatdirect}\\
		\hnsy(\hbv) &:= [\bigl(\bigl(\hnu_i\hpsi_l, N(\hat v) \bigr) \bigr)]_{i=1,\dotsc, \hat s;~l=1,\dotsc,\hat q}, \label{eq:tspodgennonl}\\
		\intertext{and}
		\hfsy &:= [\bigl(\bigl(\hnu_i\hpsi_l, f \bigr) \bigr)]_{i=1,\dotsc, \hat s;~l=1,\dotsc,\hat q}, \label{eq:gennonlin}
	\end{align}
\end{subequations}
are the Galerkin projections of the system operators and the source term assembled in the corresponding inner products.

\begin{rem}
	In the space-time Galerkin POD context, the reduced bases are projections of standard finite element bases. Concretely, by virtue of Lemma \ref{lem:optspacebas} and Lemma \ref{lem:opttimbas} one has that
	\begin{equation*}
		\hPsi = U_{\hat s}^\trp	\Ls^{-1} \Psi \quad\text{and}\quad \hUpsilon = V_{\hat q}^\trp \Ly^{-1}\Upsilon,
	\end{equation*}
where the columns of $U_{\hat s}$ and $V_{\hat q}$ are orthonormal and where $\Ls$ and $\Ly$ are factors of the mass matrices associated with $\Psi$ and $\Upsilon$. Accordingly the coefficients in \eqref{eq:spacetimecoeffs} are given as
	\begin{subequations} \label{eq:spacetimecoeffsasproj}
		\begin{align}
			\hms &:= U_{\hat s}^\trp	\Ls^{-1} \bigl[\tint \Psi \Psi^\trp \inva s\bigr] \Ls^{-\trp} U _{\hat s} = U_{\hat s}^\trp	\Ls^{-1} M_{\Oust} \Ls^{-\trp} U _{\hat s} = I_{\hat s},\\
			\hdms &:= U_{\hat s}^\trp	\Ls^{-1} \bigl[\tint \Psi \dot \Psi^\trp \inva s\bigr] \Ls^{-\trp} U _{\hat s}\label{eq:timestiffmatprj},\\
			\hmy &:= V_{\hat q}^\trp \Ly^{-1} \bigl[\sint \Upsilon \Upsilon^\trp \inva x\bigr] \Ly^{-\trp}V_{\hat q} =V_{\hat q}^\trp \Ly^{-1} M_{\Ouss} \Ly^{-\trp}V_{\hat q} = I_{\hat q}, \\
			\hay &:= V_{\hat q}^\trp \Ly^{-1} \bigl[\sint \nabla\Upsilon \nabla\Upsilon^\trp \inva x\bigr] \Ly^{-\trp}V_{\hat q}.\label{eq:spacstiffmatprj}
		\end{align}
	\end{subequations}

	Note that, despite their larger size, stiffness matrices of the standard finite element discretization, as they appear in \eqref{eq:timestiffmatprj} and \eqref{eq:spacstiffmatprj}, may be assembled much faster than the stiffness matrices $\hdms$ and $\hay$ in the formulation given in \eqref{eq:timestiffmatdirect} and \eqref{eq:spacstiffmatdirect}.
\end{rem}

%% file: chap-implementationissues.tex
\section{Implementation Issues}
In this section, we adress how to compute the measurement matrices by means of standard tools, how to incorporate the initial and terminal values in the time discretization, and how to preassemble quadratic nonlinearities.

\subsection{Computation of the Measurements}\label{sec:compthemeas}
We explain how the measurements (cf. Remark \ref{rem:measuremeasure}) that are needed for the computation of the optimal low-rank bases (cf. Lemma \ref{lem:optspacebas} and Lemma \ref{lem:opttimbas}) can be obtained in practical cases. 

In the standard \emph{method-of-lines} approach, a $\Ouss$ will be used as the FE space for a Galerkin spatial discretization that approximates \eqref{eq:genPDE} by an ODE. In a second step, a time integration scheme is employed to approximate coefficients $v_1, \dotsc, v_q \colon (0,T]\to \mathbb R^{}$ of  the solution 
\begin{equation*}
	\bar v\colon (0,T] \to \Ouss\colon t \mapsto \sum_{i=1}^qv_i(t)\nu_i
\end{equation*}
of the resulted ODE. With this and with a chosen time measurement space $\Oust$, a numerical computed measurement in $\Ousts$ of the actual solution $v$ of \eqref{eq:genPDE}, is given as 
\begin{equation}
	\mathbf{X} = 
	\begin{bmatrix} 
		(v_1, \psi_1)_\mathcal{S} & \hdots & (v_1, \psi_{\nsnap})_\mathcal{S} \\ \vdots & \ddots & \vdots \\ (v_q , \psi_1)_\mathcal{S} & \hdots & (v_q ,\psi_{\nsnap})_\mathcal{S} 
	\end{bmatrix} \Ms^{-1}.
 \label{eq:newSSmatr}
\end{equation}
% \begin{equation*}
% 	\Ms = \bigl[ (\psi_i, \psi_j)_\Oust \bigr]_{i,j=1,\dotsc,s}
% \end{equation*}

\begin{rem}
	Since in \eqref{eq:newSSmatr}, the matrix $\mathbf X$ is computed from a function with values in $\Ouss$, the inner products in $\Ouss$ and the inverse of $\My$ realizing the $L^2$-projection onto $\Ouss$ in \eqref{eq:somprojmatrix} are not present.
\end{rem}

\begin{rem}
	For smooth trajectories and for measurements using delta distributions centered at some $t_j\in (0,T)$, $j=1,\dotsc,s$, with $\tint v_i\delta(t_j)\inva t = v_i(t_j)$ the matrix \eqref{eq:newSSmatr} degenerates to the standard POD \textit{snapshot matrix}. In this case, since the delta distributions are not element of $L^2( 0,T)$, there is no way to define an optimal time basis as in Lemma \ref{lem:opttimbas}. However, one can define an optimal low-rank spatial basis by Lemma \ref{lem:optspacebas} which reduces to the standard POD optimality result with $\Ms=I$, cf. Remark \ref{rem:podisdegengenpod} and \ref{rem:genpodonlyinspace}.
\end{rem}

\subsection{Treatment of the Initial Value}\label{sec:treattheini}
The initial value \eqref{eq:genPDE_init} requires a special consideration. Firstly, like the solution of the PDE is only well defined when the initial condition is specified, also the space-time Galerkin discretized system \eqref{eq:spatimgalprojsys} is uniquely solvable if an initial condition is provided. Secondly, in particular in view of optimal control, the initial value can be subject to changes which should be realizable in the discretized model.

To maintain the prominent role of the initial condition also in the time discretization, we proceed as follows: 
\begin{enumerate}
	\item We choose an $\Oust$ that is spanned by a nodal basis $\{\psi_1,\dotsc,\psi_s\}$ and that $\psi_1$ is the basis function associated with the node at $t=0$. 
	\item For a given function, we compute $\mathbf X_0$ as in \eqref{eq:somprojmatrix} or \eqref{eq:newSSmatr} setting $\psi_1=0$ and $U_{\hat s,0}$ as the matrix of the $\hat s -1$ leading right singular vectors of $\Ly^\trp \mathbf X_0 \Ls$.
	\item We set 
		\begin{equation*}
			U_{\hat s} = \left[ \Ls^T 
			\begin{bmatrix} 
				1 \\ 0 \\ \vdots \\ 0
			\end{bmatrix} ~ U_{\hat s, 0} \right]
		\end{equation*}
		and compute the reduced basis as in Lemma \ref{lem:opttimbas} as
	\begin{equation*}
\begin{bmatrix} \hat \psi_1 \\ \hat \psi_2 \\ \vdots \\ \hat \psi_{\hat s} \end{bmatrix} 
	= U_{\hat s}^\trp \Ls ^{-1}
\begin{bmatrix}  \psi_1 \\  \psi_2 \\ \vdots \\  \psi_{s} \end{bmatrix}.
	\end{equation*}
\end{enumerate}
By this construction we obtain that $\hat \psi_1 = \psi_1$ will be associated with the initial value, whereas $\hpsi_2(0)= \dotsc= \hpsi_{\hat s}(0)=0$ will still optimally approximate the trajectory.

\subsection{Assembling of Quadratic Nonlinearities}
As an example, we consider the nonlinearity in the Burgers' equation
\begin{align} \label{eq:burgersnonlin}
 % \partial_t \brgv(t,x) + 
 \frac12 \partial_x \brgv(t,x)^2 
 % - \nu \partial_{x x} \brgv(t,x) = 0,
\end{align}
with the spatial coordinate $x\in(0,1)$, and the time variable $t\in (0,1]$.
% , and zero Dirichlet boundary conditions and an initial condition.

In the time-space Galerkin projection \eqref{eq:tspodsolansatz}, the $il$-component of the discretized nonlinearity \eqref{eq:tspodgennonl} in the case of \eqref{eq:burgersnonlin}, is given as

\begin{align*}
	H_{il}\left(\hbv\right) &= \frac 12 \int_0^1\int_0^1 \hnu_i\hpsi_l\cdot \partial_x \hat v^2 \inva x \inva t \\
	&= \frac 12 \int_0^1\int_0^1 \hnu_i\hpsi_l\cdot \partial_x(([\hsbas^\trp \otimes \hybas^\trp] \hbv)^2) \inva x \inva t  \\
	&= \hbv^\trp[\int_0^1\hnu_i \hsbas \hsbas^\trp \inva t \otimes \frac 12 \int_0^1 \hpsi_l\partial_x(\hybas \hybas^\trp)^2\inva x] \hbv,
	% &= \hbv^\trp[\int_0^1\hnu_i \hsbas \hsbas^\trp \inva t \otimes \frac 12 \int_0^1 \hpsi_l(\hybas\partial_x \hybas^\trp+ \partial_x(\hybas) \hybas^\trp)\inva x] \hbv
\end{align*}
where we have used the linearity of the \emph{Kronecker product} and that 
\begin{equation*}
	\hat v^2 = ([\hsbas^\trp \otimes \hybas^\trp] \hbv)^2 = \hbv^\trp [\hsbas \otimes \hybas][\hsbas^\trp \otimes \hybas^\trp]\hbv = \hbv^\trp [\hsbas\hsbas^\trp \otimes \hybas \hybas ^\trp] \hbv.
\end{equation*}
Thus, the evaluation of the discretized nonlinear term can be assisted by precomputing 
\begin{equation*}
	\int_0^1\hnu_i \hsbas \hsbas^\trp \inva t \andi \frac 12 \int_0^1 \hpsi_l(\hybas\partial_x \hybas^\trp+ \partial_x(\hybas) \hybas^\trp)\inva x
\end{equation*}
for all $\hnu_i$, $i=1, \dotsc, \hat s$ and $\hpsi_l$, $l=1, \dotsc , \hat q$.
\begin{rem}
	If $V_{\hat q}$ is the matrix of the spatial POD modes that transform the FEM basis $\ybas$ into the reduced basis $\hybas$ via $\hybas = V_{\hat q}^\trp\Ly^{-1} \ybas$, then the spatial part of the reduced nonlinearity fulfills
\begin{align*}
	\frac 12 \int_0^1 \hpsi_l(\hybas\partial_x \hybas^\trp+ \partial_x(\hybas) \hybas^\trp)\inva x = \\
	\frac 12 V_{\hat q}^\trp \Ly^{-1}\int_0^1 \hpsi_l(\ybas\partial_x \ybas^\trp+ \partial_x(\ybas) \ybas^\trp)\inva x \Ly^{-\trp}V_{\hat q},
\end{align*}
where the inner matrix of the latter expression might be efficiently assembled in a FEM package. The same idea applies to the time-related part.

\end{rem}

%% file: chap-optcontrol.tex
% \section{Space-time POD for Optimal Control}
\section{Application in PDE-Constrained Optimization}
We consider a generic optimal control problem.
\begin{prob}\label{prob:genoptcont}
For a given target trajectory $x^* \in L^2( 0,T;L^2( \Omega))$ and a penalization parameter $\alpha>0$, we consider the optimization problem
\begin{equation}
	\mathcal J(x,u):=\frac 12 \norm{x-x^*}_{L^2}^2+\frac \alpha2 \norm{u}_{L^2}^2 \to \min_{u\in L^2(0,T;L^2( \Omega))}
	\label{eq:costfunc}
\end{equation}
subject to the generic PDE
\begin{subequations}\label{eq:genPDEcontr}
	\begin{align}
		\dot x - \Delta x + N(x) &= f + u\quad \text{on } (0,T] \times \Omega, \\
		x\bigr|_{\partial \Omega} &= 0 \quad \hspace{0.63cm}\text{on } (0,T],\\
		x\bigr|_{t=0} &= x_0 \quad \hspace{0.43cm}\text{on } \Omega.
	\end{align}
\end{subequations}
\end{prob}

If the nonlinearity is smooth, then necessary optimality conditions with respect to Problem \ref{prob:genoptcont} for $(x,u)$ are given through $u=\frac 1\alpha \lambda$, where $\lambda$ solves the adjoint equation
\begin{subequations}\label{eq:genPDEadjcontr}
	\begin{align}
		-\dot \lambda - \Delta \lambda + D_xN(x)^\trp \lambda + x&= x^* \quad \text{on } (0,T] \times \Omega, \\
		\lambda\bigr|_{\partial \Omega} &= 0 \quad \hspace{0.2cm}\text{on } (0,T],\\
		\lambda\bigr|_{t=T} &= 0 \quad \hspace{0.22cm} \text{on } \Omega,
	\end{align}
\end{subequations}
where $D_x$ denotes the Frech\'et derivative, which is coupled to the state equation \eqref{eq:genPDEcontr} through $x$ and $u$; see \cite{Tr09}.

Given low-dimensional spaces $\hOust:=\spann\{\hat \psi_1, \dotsc, \hat \psi_{\hat s}\}$, $\hAdjst := \spann\{\hphi_1, \dotsc, \hphi_{\hat r}\} \subset H^1( 0,T)$ and $\hOuss:=\spann\{\hat \nu_1, \dotsc, \hat \nu_{\hat q}\}$, $\hAdjss := \spann\{\lambda_1, \dotsc, \lambda_{\hat p}\}\subset H_0^1( \Omega)$, a tensor space-time Galerkin discretization of the coupled system \eqref{eq:genPDEcontr}-\eqref{eq:genPDEadjcontr} reads
\begin{subequations}\label{eq:spatimgalprojoptsys}
\begin{align} 
[\hdms \otimes \hmy +\hms \otimes \hay ]\hbv + \hnsy(\hbv) -\frac 1\alpha[\hmsr \otimes \hmyl] \hadjvd &= \hfsy,\label{eq:spatimgalprojoptsys_state}\\
[-\hdmr \otimes \hml +\hmr \otimes \hal ]\hadjvd + \hnrl(\hbv)\hadjvd +[\hmrs \otimes \hmly] \hbv &= [\hmrs \otimes \hmly] \hbv^*,
\end{align}
\end{subequations}
with the coefficients $\hdmr$, $\hmr$, $\hml$, $\hal$ and the nonlinearity $\hnrl(\hbv)\hadjvd$ defined as in \eqref{eq:spatimgalprojsys}, with $\hmsr$, $\hmrs$, $\hmyl$, $\hmly$ denoting the mixed mass matrices like
\begin{equation*}
	\hmsr := [\bigl(\hpsi_\ell, \hphi_k \bigr)]^{\ell = 1, \dotsc, \hat s}_{k=1, \dotsc, \hat r} \in \mathbb R^{\hat s, \hat r} ,
\end{equation*}
with $\hbv^*$ representing the target $v^*$ projected onto $\hOusts$, with the spatial boundary conditions resolved in the ansatz spaces, and with accounting for the initial and terminal conditions via requiring
\begin{equation*}
	\hat v(0) = \sum_{j=1}^{\hat s} \sum_{i=1}^{\hat q} \xtsvecij^1 \hat \nu_i \hat \psi_j(0) = \Pi_{\hOuss}x_0 \quad\text{and}\quad \hat \lambda(T) = \sum_{j=1}^{\hat r} \sum_{i=1}^{\hat p} \xtsvecij^1 \hmu_i \hphi_j(T) = 0,
\end{equation*}
cf. Chapter \ref{sec:treattheini}.

% \begin{table}[ht]
%  \centering
%  \begin{tabular}{l|c|c|c}
%   variable & space & bases & dimension \\
%   \hline
%   state $x$ & $\mathcal{Y} \cdot \mathcal{S}$ & $\{\nu_1,...,\nu_q\} \times \{\psi_1,...,\psi_s\}$ & $qs$\\
%   reduced state $\hat x$ & $\hat{\mathcal{Y}} \cdot \hat{\mathcal{S}}$ & $\{\hat \nu_1,...,\hat\nu_{\hat q}\} \times \{\hat\psi_1,...,\hat\psi_{\hat s}\}$ & $\hat q \hat s$\\
%   control $u$, adjoint $\lambda$ & $\Lambda \cdot \mathcal{R}$ & $\{\lambda_1,...,\lambda_p\} \times \{\phi_1,...,\phi_r\}$ & $pr$\\
%   reduced control/adjoint $\hat u, \hat \lambda$ & $\hat \Lambda \cdot \hat{\mathcal{R}}$ & $\{\hat{\lambda}_1,...,\hat{\lambda}_{\hat{p}}\} \times \{\hat{\phi}_1,...,\hat{\phi}_{\hat{r}}\}$ & $\hat{p}\hat{r}$
%  \end{tabular}
%  \caption{Overview on the used function spaces, with $\hat q \hat s \ll qs$ and $\hat p \hat r \ll pr$.} \label{tab:used_spaces}
% \end{table}

%% file: chap-numerics.tex
\section{Numerical Experiments}
We consider the optimal control of a Burgers' equation as it was described in \cite{H08,KV99}. Therefore, in Problem \ref{prob:genoptcont}, we replace the generic PDE \eqref{eq:genPDEcontr} by Burgers' equation, namely:
\begin{subequations}\label{eq:burgeroptcont}
	\begin{align}
	  \dot x - \nu\partial_{\xi\xi} x + \frac12 \partial_\xi (x^2) &= u \hspace{0.52cm} \text{on } (0,T] \times (0, L), \\
	  x\bigr|_{\xi=0,\xi=L} &= 0  \hspace{0.53cm} \text{on } (0,T],\\
		x\bigr|_{t=0} &= x_0 \quad \text{on } (0,L),
	\end{align}
\end{subequations}
where $L$ and $T$ denote the length of the space and time interval and where $\nu > 0$ is a parameter. We set $T=1$ and $L=1$ and, as the initial value, we take the step function 
\begin{equation}\label{eq:inival}
  x_0\colon (0,1) \to \mathbb R^{}\colon \xi \mapsto
  \begin{cases}
	1, \quad \text{if }\xi \leq 0.5 \\
	0, \quad \text{if }\xi < 0.5
  \end{cases}.
\end{equation}
As for the target, we define $x^*$ via $x^*(t)=x_0$, which means that the optimization is designed to keep the system in its initial state.

Thus the concrete problem reads as follows:
\begin{prob}
\label{prob:optcon}
  Given parameters $\nu$ and $\alpha$, find $u\in L^2(0,1; L^2( 0,1))$ such that 
\begin{equation}\label{eq:burgercostfunc} 
	\frac 12 \int_0^1 \int_0^1 (x(t, \xi)-x_0)^2\inva \xi \inva t+\frac \alpha 2 \int_0^1 \int_0^1 u^2(t, \xi)\inva \xi \inva t \to \min_{u\in L^2(0,1;L^2( 0,1))}
\end{equation}
subject to Burgers' equation \eqref{eq:burgeroptcont}.
\end{prob}
\subsection{Space-time Generalized POD for Optimal Control} \label{sec_spatimpod_burgers}
The general procedure is as follows:
\begin{enumerate}
  \item Do at least one forward solve of the state equation \eqref{eq:burgeroptcont} and at least one backward solve of the adjoint equation \eqref{eq:genPDEadjcontr} to setup generalized measurement matrices of the state and the costate as explained in Section \ref{sec:compthemeas}.
  \item Compute optimized space and time bases for the state and the costate as defined in Lemma \ref{lem:optspacebas} and \ref{lem:opttimbas}. To account for the initial and the terminal value, one may resort to the procedure explained in Section \ref{sec:treattheini}.
  \item Set up the projected closed-loop optimality system \eqref{eq:spatimgalprojoptsys} and solve for the optimal costate $\hadjvd$ of the reduced system.
  \item Lift $\hat{\boldsymbol u} = \frac 1\alpha \hadjvd$ up to the full space-time grid and apply it as suboptimal control to the actual problem.
\end{enumerate}

The procedure is defined by several parameters. In the presented examples, we fix $\mathcal Y = \Lambda $ and $\mathcal S = \mathcal R$, corresponding to the initial space and time discretizations, and investigate the influence of the other parameters on the numerical solution of the optimal control problem. See Table \ref{tab_genpod_parameters} for an overview of the parameters and their default values.

\begin{table}[tb]
  \begin{tabular}{|c|p{.4\textwidth}|c|c|}
	\hline
  Parameter & Description & Base Value & Range \\
  \hline
  $\mathcal Y$, $\Lambda $ & Space of piecewise linear finite elements on an equidistant grid of dimension $q$, $p$ & $q =p  = 220$ & -- \\
  $\mathcal S$, $\mathcal R $ & Space of linear hat functions on an equidistant grid of dimension $s$, $r$ & $s =r  = 120$ & -- \\
  $\hOuss$, $\hat \Lambda $ & POD reductions of $\mathcal Y$ and $\Lambda$ of dimension $\hat q$, $\hat p$; cf. Lemma \ref{lem:optspacebas} & $\hat q = \hat p  = 12 $ & $6$ -- $24$ \\
  $\hOust$, $\hAdjst$ & POD reductions of $\mathcal S$ and $\mathcal R$ of dimension $\hat s$, $\hat r$; cf. Lemma \ref{lem:opttimbas} & $\hat r = \hat s  = 12 $ & $6$ -- $24$ \\
  $\alpha$ & Regularization parameter in the cost functional \eqref{eq:burgercostfunc}& $1\cdot 10^{-3}$ & $2.5\cdot 10^{-4}$ -- $1.6\cdot 10^{-2}$  \\
  $\nu$ & Viscosity parameter in the PDE  \eqref{eq:burgeroptcont}& $2\cdot 10^{-3}$ & $5\cdot 10^{-4}$ -- $1.6\cdot 10^{-2}$ \\
  \hline
 \end{tabular}
 \caption{Description and values of the parameters of the numerical examples of Section \ref{sec_spatimpod_burgers}}
  \label{tab_genpod_parameters}
\end{table}

We will measure the performance of the approach through:
\begin{itemize}
  \item The time $\ctime$ it takes to solve the reduced optimality system \eqref{eq:spatimgalprojoptsys} for $\hadjvd$, reporting the best number out of $5$ runs.
  \item The difference $\norm{\hat x - x_0}_{L^2}$ between the target state and the state $\hat x$ achieved by using the suboptimal control $\hat u$ on the base of $\hadjvd$ in the simulation of the full model.
  \item The value $\mathcal J(\hat x, \hat u)$ of the cost function \eqref{eq:burgercostfunc}.
\end{itemize}
The spatial discretization is carried out with the help of the FEM library \emph{FEniCS} \cite{LogOeRW12}. For the time integration, we use \emph{SciPy}'s builtin ODE-integrator \texttt{scipy.integrate.odeint}. To solve the nonlinear system \eqref{eq:spatimgalprojoptsys} for $\hadjvd$, we use \emph{SciPy}'s routine \texttt{scipy.optimize.fsolve}. The norms are approximated in the used FEM space. 
%The reported timings \ctime~are the best out of 5 consecutive runs. 
The implementation and the code for all tests as well as the documentation of the hardware are available from the author's public git repository \cite{Hei15_gh-spacetime-genpod-burgers}; see also the section on code availability on page \pageref{sec:codeavailability}.

{\bf Choice of the measurements}. The computation of the measurements and the choice of the reduced bases is an important parameter of the approach.
  Generally the basis of $\hOust \cdot \hOuss$ should be well suited to approximate the state, whereas the basis $\hAdjst\cdot \hAdjss$ should well represent the adjoint state. In the optimization case, where the suboptimal input is defined through $\frac 1\alpha \hadjvd$ and its lifting to the full-order space, two other conditions emerge. Firstly, the reduced basis of the adjoint state, should also well approximate the optimal control. Secondly, the bases of the state and the adjoint must not be orthogonal or ``almost'' orthogonal such that the joint mass matrix $[\hmsr \otimes \hmyl]$ degenerates and the contribution of the input in \eqref{eq:spatimgalprojoptsys_state} vanishes.

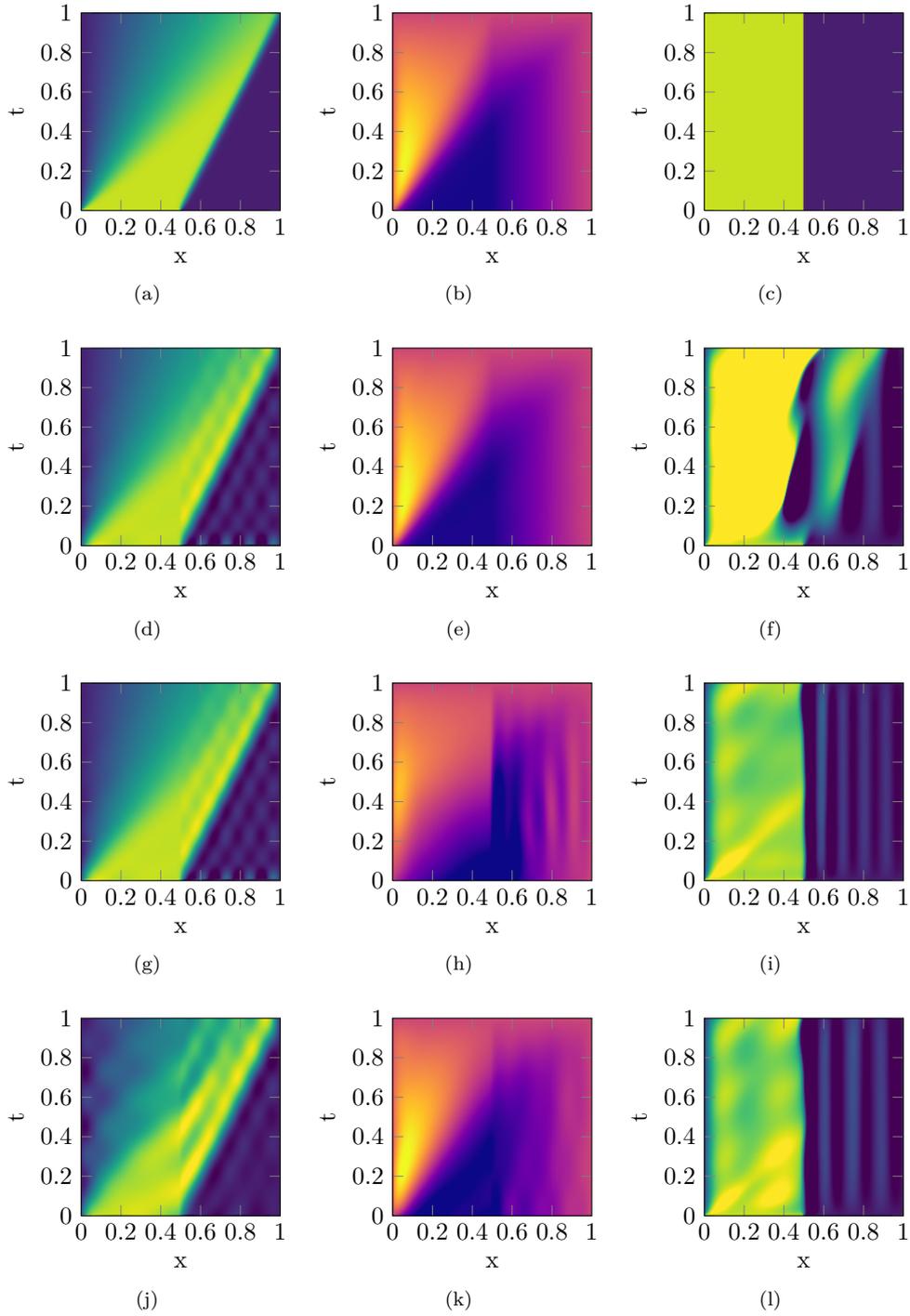
\begin{figure}[p]
	\begin{center}
	  \subfigure[][]{\input{plots/1.tikz}}
	  \subfigure[][]{\input{plots/2.tikz}}
	  \subfigure[][]{\input{plots/3.tikz}}
	  \subfigure[][]{\input{plots/4.tikz}}
	  \subfigure[][]{\input{plots/5.tikz}}
	  \subfigure[][]{\input{plots/6.tikz}}
	  \subfigure[][]{\input{plots/7.tikz}}
	  \subfigure[][]{\input{plots/8.tikz}}
	  \subfigure[][]{\input{plots/9.tikz}}
	  \subfigure[][]{\input{plots/10.tikz}}
	  \subfigure[][]{\input{plots/11.tikz}}
	  \subfigure[][]{\input{plots/12.tikz}}
	\end{center}
	\caption{Illustration of the effect of the choice of the snapshots on the performance of the low-dimensional approximations of the state (a), the adjoint state (b), and on how outcome of the optimization matches the target state (c). The second row (d-f) corresponds to the case that snapshots of the state and the adjoint are used to approximate the state and the adjoint, respectively. For the results depicted in the third row (g-i), the optimized basis for the state was used also for the adjoint. The results depicted in the last row (j-l) were obtained by combining state and adjoint snapshots for the computation of the reduced bases. For a comparable illustration we have used color maps with linear intensity on the intervals $[-0.1, 1.1]$ for the states and $[-0.5, 0.5]$ for the adjoint states. Values that exceeded these margins were cropped.}
	\label{fig:allNCG}
\end{figure}

  As illustrated in the plots in Figure \ref{fig:allNCG}, the straight-forward approach of constructing the bases for the state by means of state measurements and the basis for the adjoint by means of measurements of the adjoint, well approximates the state and the adjoint but not the closed-loop problem. It turned out that taking the state measurements to also construct the reduced space for the adjoint gave a better approximation to the optimality system while, naturally, only poorly approximating the adjoint. The best result were obtained in combining state and adjoint state measurements to construct the bases. 

Thus, for the computation of the optimal bases for the following tests, we combined the measurements obtained from one forward solve with no control and one backward solve with the state from the forward solve and the target state. 

\newlength\figureheight
\newlength\figurewidth

\setlength\figureheight{4.4cm}
\setlength\figurewidth{4.4cm}

{\bf Overall Number of Modes.}
In these tests we examine how the number of modes, i.e. the dimension of the reduced order system, affects the time needed to solve the reduced order system and the value of the cost functional \eqref{eq:burgercostfunc} achieved through the suboptimal control. 

We fix $\nu=0.005$ and set $\hatq=\hatp=\hatr=\hats=\hat K/4$, with $\hat K \in \{24, 36, 48, 72, 96\}$. Thus, for every setup, the nonlinear system \eqref{eq:spatimgalprojoptsys} of dimension $\hat K$ has to be solved for the optimal costate $\hadjvd$. The results of these tests are reported in Table \ref{tab:numpodmodes}.

\begin{table}[p]
  \begin{tabular}{l|ccccc}
$\hat K$    &  $24$ & $36$ & $48$ & $72$ & $96$ \\
\hline
$\vterm $                       &  $0.0330$ & $0.0280$ & $0.0192$ & $0.0121$ & $0.0104$ \\
$\mathcal J(\hat x, \hat u)$    &  $0.0351$ & $0.0309$ & $0.0234$ & $0.0177$ & $0.0152$ \\
\ctime~$[s]$    & $\texttt{0.1}$ & $\texttt{0.48}$ & $\texttt{1.81}$ & $\texttt{18.7}$ & $\texttt{155}$
% $0.1624$ & $0.7762$ & $2.9026$ & $29.3518$ & $223.5538$
  \end{tabular}
  \caption{Performance of the suboptimal control versus varying resolutions of space and time.}
  \label{tab:numpodmodes}
\end{table}

As expected, the larger the reduced model, the lower the achieved values of the cost functional. Also, with growing order of the reduced model, the time needed to solve the corresponding nonlinear system increases drastically. %In the considered setup this increase is drastic; when doubling the number of modes, then computation time grows by a factor of $100$.

{\bf Space vs. Time Reduction} From the previous tests, we found that in the considered setup, an overall number of $\hat K=48$ modes is a good compromise between accuracy and computation time. In this section, we examine how the distribution of modes between space and time affects the quality of the suboptimal control. Therefore, and for a varying increment/decrement $j$, we set $\hat q = \hat p := 12\mp j$ and $\hat s =\hat r := 12 \pm j$. Accordingly, the overall number of degrees of freedom stays constant throughout the tests but we add weight on the approximation of either the time or the space component. 

The results are listed in Table \ref{tab:numpodmodedist}. It turns out that up to a certain level, it is beneficial to emphasize on the space as opposed to the time resolution. In the considered setup, the distribution of $16$ degrees of freedom in space versus $8$ in time gave the best results in terms of performance of the corresponding suboptimal control. Interestingly, also the timings \ctime~vary significantly. This variance is due to different convergence behavior of the optimization algorithm used to solve the nonlinear system.

\begin{table}[p]
  \begin{tabular}{l|ccccccc}
$(\hat q, \hat s)$/$(\hat p, \hat r)$ & (18,  6) &(17,  7) &(16,  8) & (14, 10) & (12, 12) & (10, 14) & ( 8, 16) \\
\hline
$\vterm$                     & $0.0138$ & $0.0125$ & $0.0117$ & $0.0137$ & $0.0192$ & $0.0326$ & $0.0339$ \\
$\mathcal J(\hat x, \hat u)$ & $0.0184$ & $0.0173$ & $0.0167$ & $0.0184$ & $0.0234$ & $0.0364$ & $0.0364$ \\
\ctime & $\texttt{0.98}$ & $\texttt{1.84}$ & $\texttt{1.19 }$ & $\texttt{1.57}$ & $\texttt{1.81}$ & $\texttt{1.49 }$ & $\texttt{1}$
  \end{tabular}
  \caption{Performance of the suboptimal control versus varying distributions of space and time resolutions.}
  \label{tab:numpodmodedist}
\end{table}

{\bf Reduced Order Model vs. Viscosity Parameter} In these tests, we examine how the low-rank space-time Galerkin approach performs over a range of viscosity parameters $\nu$. It is known that for low values of $\nu$, the problem is \emph{convection dominated} and hard to approximate by POD bases. In the considered setup, where the target state is a nonsmooth function, we also expect a decreasing performance for larger values of $\nu$, since the diffusion makes the step in the target untrackable.

The results are listed in Table \ref{tab:nuperform168} for $(\hat q, \hat s)=(\hat p, \hat r)=(16,8)$, which was the optimal distribution as found in the previous tests, and in Table \ref{tab:nuperform1212} for $(\hat q, \hat s)=(\hat p, \hat r)=(12,12)$. The optimal distribution $(16,8)$ has its performance peak at $\nu = 8\cdot 10^{-3}$. It also shows the expected pattern that for lower values of $\nu$, for which the nonlinearity gets more emphasis, the time for the solution of the nonlinear system increases. The results of the runs with equally distributed numbers of space and time modes $(12,12)$ are listed in Table \ref{tab:nuperform1212}. The comparison to the distribution $(16,8)$ does not give a clear conclusion. For higher viscosities the $(12, 12)$-case is worse in all categories. For low viscosities, it outperforms the $(16,8)$ case in terms of computation time and, for the lowest investigated value of $\nu$, even in terms of approximation quality.

\begin{table}[p]
  \begin{tabular}{l|ccccccc}
$\nu$   & $5\cdot 10^{-4}$ & $1\cdot 10^{-3}$ & $2\cdot 10^{-3}$ & $4\cdot 10^{-3}$ & $8\cdot 10^{-3}$ & $1.6\cdot 10^{-2}$ & $3.2\cdot 10^{-2}$ \\
\hline
$\vterm $  & $1.1072$ & $0.1368$ & $0.0188$ & $0.0126$ & $0.0098$ & $0.0111$ & $0.0198$\\
$\cfune$   & $1.8489$ & $0.1994$ & $0.0225$ & $0.0173$ & $0.0150$ & $0.0168$ & $0.0268$\\
\ctime & $\texttt{2.95}$ & $\texttt{7.13}$ & $\texttt{1.91}$ & $\texttt{1.04}$ & $\texttt{1.07}$ & $\texttt{1}$ & $\texttt{1.64}$

  \end{tabular}
  \caption{Performance of the suboptimal control versus varying diffusion parameters $\nu$ for $(\hat q, \hat s)=(\hat p, \hat r)=(16,8)$.}
  \label{tab:nuperform168}
\end{table}

\begin{table}[p]
  \begin{tabular}{l|ccccccc}
$\nu$   & $5\cdot 10^{-4}$ & $1\cdot 10^{-3}$ & $2\cdot 10^{-3}$ & $4\cdot 10^{-3}$ & $8\cdot 10^{-3}$ & $1.6\cdot 10^{-2}$ & $3.2\cdot 10^{-2}$ \\
\hline
\vterm & $0.0236$ & $0.0259$ & $0.0263$ & $0.0217$ & $0.0153$ & $0.0123$ & $0.0210$\\
\cfune & $0.0269$ & $0.0293$ & $0.0299$ & $0.0256$ & $0.0201$ & $0.0176$ & $0.0281$\\
\ctime	& $\texttt{1.60}$ & $\texttt{2.06}$ & $\texttt{1.76}$ & $\texttt{1.24}$ & $\texttt{1.64}$ & $\texttt{1.8}$ & $\texttt{1.9}$
  \end{tabular}
  \caption{Performance of the suboptimal control versus varying diffusion parameters $\nu$ for $(\hat q, \hat s)=(\hat p, \hat r)=(12,12)$.}
  \label{tab:nuperform1212}
\end{table}

{\bf Regularization Parameter} In this section, we examine the influence of the regularization parameter $\alpha$ onto the performance. Therefore we fix $\nu=5\cdot10^{-3}$, $(\hat q, \hat s)=(\hat p, \hat r)\in\{(16,8), (12, 12)\}$, and vary $\alpha$ which defines the penalization of the control magnitude in the cost functional \eqref{eq:burgercostfunc}. 

The results are listed in Table \ref{tab:alpha168} for the $(16,8)$-case and in Table \ref{tab:alpha1212} for the $(12,12)$-case. For higher values of $\alpha$, for both cases, the optimization performs similarly well in approximating the control problem, with a clear advantage of the $(16,8)$ distribution in terms of computation time. For lower values of $\alpha$ the $(16,8)$ setup outperforms the $(12, 12)$ by far. Interestingly, for the smallest $\alpha$, the tracking error \vterm~increases again.

\begin{table}[p]
  \begin{tabular}{l|ccccccc}
$\alpha$ & $2.5\cdot 10^{-4}$ & $5\cdot 10^{-4}$ & $1\cdot 10^{-3}$ & $2\cdot 10^{-3}$ & $4\cdot 10^{-3}$ & $8\cdot 10^{-3}$ & $1.6\cdot 10^{-2}$\\
\hline
\vterm & $0.0137$ & $0.0121$ & $0.0117$ & $0.0124$ & $0.0144$ & $0.0179$ & $0.0237$\\
\cfune & $0.0158$ & $0.0155$ & $0.0167$ & $0.0196$ & $0.0240$ & $0.0305$ & $0.0398$\\
\ctime & $\texttt{1.12}$ & $\texttt{0.96}$ & $\texttt{1.18}$ & $\texttt{1.4}$ & $\texttt{1.59}$ & $\texttt{1.82}$ & $\texttt{1.88}$
  \end{tabular}
  \caption{Performance of the suboptimal control versus varying regularization parameters $\alpha$ for $(\hat q, \hat s)=(\hat p, \hat r)=(16,8)$.}
  \label{tab:alpha168}
\end{table}

\begin{table}[p]
  \begin{tabular}{l|ccccccc}
$\alpha$ & $2.5\cdot 10^{-4}$ & $5\cdot 10^{-4}$ & $1\cdot 10^{-3}$ & $2\cdot 10^{-3}$ & $4\cdot 10^{-3}$ & $8\cdot 10^{-3}$ & $1.6\cdot 10^{-2}$\\
\hline
\vterm & $0.0488$ & $0.0293$ & $0.0192$ & $0.0148$ & $0.0145$ & $0.0168$ & $0.0215$\\
\cfune & $0.0504$ & $0.0318$ & $0.0234$ & $0.0213$ & $0.0239$ & $0.0301$ & $0.0393$\\
\ctime & $\texttt{1.31}$ & $\texttt{1.59}$ & $\texttt{1.81}$ & $\texttt{1.97}$ & $\texttt{2.02}$ & $\texttt{2.5}$ & $\texttt{3.73}$
  \end{tabular}
  \caption{Performance of the suboptimal control versus varying regularization parameters $\alpha$ for $(\hat q, \hat s)=(\hat p, \hat r)=(12,12)$.}
  \label{tab:alpha1212}
\end{table}
\subsection{Gradient-based Classical POD-reduced Optimal Control}
In this section we consider the suboptimal numerical solution of Problem \ref{prob:optcon} based on a (spatial) POD reduction in a classical \textit{method-of-lines} approach, cf. \cite{KV99}. After finite element discretization in space, we consider \eqref{eq:burgeroptcont} in semi-discretized form,
% \begin{align*}
% \min_{\mathbf{u}} \int_0^T  \frac12 \mathbf{x}(\tau)^\trp \My \mathbf{x}(\tau) - ... + \frac{\alpha}{2} \mathbf{u}(\tau)^\trp \My \mathbf{u}(\tau) \inva \tau,
% \end{align*}
% subject to
\begin{subequations}\label{eq:fom_semi}
\begin{align}
 M_{\mathcal Y} \partial_t \mathbf{x}(t) + \nu K_{\mathcal Y} \mathbf{x}(t) +  H_{\mathcal Y}(\mathbf{x}) - M_{\Lambda} \mathbf{u}(t) &= \mathbf{0}, \\
 \mathbf{x}(0) &= \mathbf{x}_0,
\end{align}
\end{subequations}
where mass and stiffness matrices are defined in the same way as their reduced-order counterparts in \eqref{eq:timemassmatdirect}-\eqref{eq:gennonlin}. Classical POD is based on an SVD of the so-called snapshot matrix taken from $s$ distinct time instances,
\begin{align}
% \label{eq:snapshotmatr}
X &= [\mathbf{x}(t_1),...,\mathbf{x}(t_s)] \in \mathbb{R}^{q \times s}. \label{svd_x}%\\
% U &= [\mathbf{u}(\tau_1),...,\mathbf{u}(\tau_s)] \in \mathbb{R}^{q \times s}\label{svd_u}.
\end{align}
For $U_{\hat q}$ being the matrix that consists of the $\hat q$ leading left singular vectors of \eqref{svd_x}, we introduce the reduced state variable via $\mathbf{x}(t) \approx U_{\hat q} \hat{\mathbf{x}}(t)$. Similarly, a reduced-order control variable $\hat{\mathbf{u}}$ can be introduced. A suboptimal solution to Problem \ref{prob:optcon} can be obtained by minimizing the corresponding fully discretized POD-reduced Lagrangian function,
\begin{align}
 & \hat{\mathcal L}(\hat{\mathbf{x}}_0,...,\hat{\mathbf{x}}_{n_t},\hat{\mathbf{u}}_0,...,\hat{\mathbf{u}}_{n_t},\hat{\boldsymbol{\lambda}}_1,...,\hat{\boldsymbol{\lambda}}_{n_t}) \notag \\
  & = \sum_{j=0}^{n_t} \delta \! t \left(\frac12 \hat{\mathbf{x}}_j^\trp \hmy \hat{\mathbf{x}}_j - (\hat{\mathbf{x}}^\ast_j)^\trp \hat{\mathbf{x}}_j + \frac{\alpha}{2} \hat{\mathbf{u}}_j^\trp M_{\hat \Lambda} \hat{\mathbf{u}}_j \right) \notag \\
  & + \sum_{j=0}^{n_t-1} \hat{\boldsymbol{\lambda}}_{j+1}^\trp \left[ \hmy \left( \frac{\hat{\mathbf{x}}_{j+1} - \hat{\mathbf{x}}_j}{\delta \! t} \right) + \nu \hay \hat{\mathbf{x}}_{j+1} + H_{\mathcal Y}(U_{\hat q}\hat{\mathbf{x}}_{j+1}) - M_{\hat \Lambda} \hat{\mathbf{u}}_{j+1}\right], \label{eq:Lagr_red}
\end{align}
where we refer the reader to \cite{B13,H08} for a detailed derivation. In the following numerical experiments, we use $n_t$ time steps for an implicit Euler time integration, and use a gradient-based optimization scheme for minimizing \eqref{eq:Lagr_red}. The gradient $\nabla_{\hat u} \hat{\mathcal L}$ is computed using the adjoint approach of \cite[Algorithm 6.1]{H08} which is solved \textit{backwards in time}. To compare with the results in Section \ref{sec_spatimpod_burgers} as best as possible, we use the same base parameters, namely $s=120$ equidistantly distributed snapshots in \eqref{svd_x}, and a linear finite element basis of dimension $q = 220$ for the spatial component of the full-order model \eqref{eq:fom_semi}. In Table \ref{tab:numpodmodes}, we use the BFGS implementation \cite{Imm} with a stopping criterion that targets the objective function values achieved in Table \ref{tab:numpodmodes}. When $\hat q = n_t$ are increased, we observe that the objective funtion value improves at a linear cost. These results can be improved when DEIM \cite{DEIM} at a fixed dimension of $25$ DEIM points is used to approximation the nonlinear term in \eqref{eq:Lagr_red}. A comparison with the gradient-based SPG method \cite{MBR00} is given in \cite{B13} and yields similar results.

\begin{table}[p]
  \begin{tabular}{l|c@{\hskip 2.4mm}c@{\hskip 2.4mm}c@{\hskip 2.4mm}c@{\hskip 2.4mm}c|c@{\hskip 2.4mm}c@{\hskip 2.4mm}c@{\hskip 2.4mm}c@{\hskip 2.4mm}c}
%   \hline
  & \multicolumn{5}{ c| }{POD} & \multicolumn{5}{ c }{POD-DEIM}\\
$\hat q = n_t$    &  $6$ & $9$ & $12$ & $18$ & $24$ &  $6$ & $9$ & $12$ & $18$ & $24$\\
\hline
\vterm & \small{$0.0355$}  & \small{$0.0306$} & \small{$0.0216$} &  \small{$0.0130$} &  \small{$0.0131$} & \small{$0.0353$}  &\small{$0.0305$}  & \small{$0.0223$} &  \small{$0.0157$} &\small{$0.0103$}\\
\cfune    & \small{$0.0363$} &  \small{$0.0318$} &  \small{$0.0240$} & \small{$0.0178$} & \small{$0.0175$}& \small{$0.0356$}  & \small{$0.0310$} & \small{$0.0234$} & \small{$0.0182$} &\small{$0.0153$} \\
\texttt{\#BFGS}   &  \small{45} & \small{62} & \small{93}& \small{134} & \small{203}& \small{30}   &\small{47}  & \small{77} & \small{117} & \small{168}\\
\ctime   & \small{\texttt{0.20}} & \small{\texttt{0.41}}& \small{\texttt{0.80}} & \small{\texttt{1.78}} & \small{\texttt{3.80}}&  \small{\texttt{0.13}} &\small{\texttt{0.28}}  & \small{\texttt{0.60}} & \small{\texttt{1.43}} & \small{\texttt{2.81}}\\
% \hline
 \end{tabular}
  \caption{Performance of suboptimal control based on a classical spatial POD and POD-DEIM reduction for the state and control variable. BFGS is terminated based on targeting the suboptimal value of \cfune~ in Table \ref{tab:numpodmodes}. We fix $\nu = 0.005$ and $\alpha = 0.001$ .}
  \label{tab:pod_bfgs}
\end{table}

In a second experiment we mimic the experiments reported in Table \ref{tab:numpodmodedist}. Therefore, we vary the ratio of POD dimension $\hat p$ versus number of Euler time integration steps $n_t$. Also here, the best results are obtained when the spatial dimension is large compared to the time discretization.
\begin{table}[p]
  \begin{tabular}{l|ccccccc}
$(\hat q, n_t)$ & (18,  6) &(17,  7) &(16,  8) & (14, 10) & (12, 12) & (10, 14) & (8, 16) \\
\hline
\vterm & $0.0217$ & $0.0216$ &$0.0209$& $0.0218$& $0.0216$& $0.0205$& $0.0184$\\
\cfune &  $0.0238$ & $0.0236$ & $0.0236$& $0.0237$& $0.0240$& $0.0236$&  $0.0235$\\
\texttt{\#BFGS}   & $70$ & $72$ &$85$ & $84$ & $93$& $95$  & $111$ \\
\ctime& $\texttt{0.37}$ &  $\texttt{0.43}$ &  $\texttt{0.55}$ &$\texttt{0.66}$& $\texttt{0.80}$&$\texttt{0.88}$ & $\texttt{1.10}$
  \end{tabular}
  \caption{Performance of the POD-Lagrangian suboptimal control for varying spatial POD reduction and temporal integration.}
  \label{tab:numpodvarst}
\end{table}

In the previous two experiments, BFGS was terminated based on a priori knowledge. In Table \ref{tab:numbfgsterm}, we chose a stopping criterion based on a tolerance for the gradient of the Lagrangian \eqref{eq:Lagr_red}. The timings presented in Table \ref{tab:numbfgsterm} are for the case $\hat q = n_t = 18$ and indicate that the configuration used in Table \ref{tab:pod_bfgs} is close to optimal.
\begin{table}[p]
  \begin{tabular}{l|ccccccc}
$\texttt{tol}_\nabla$ & \texttt{1e-2} &\texttt{5e-3} &\texttt{1e-3} & \texttt{5e-4} & \texttt{1e-4} & \texttt{5e-5} & \texttt{1e-5} \\
\hline
\cfune & $0.0738$  & $0.0738$ & $0.0487$&  $0.0487$& $0.0173$ & $0.0163$ & $0.0162$ \\
\texttt{\#BFGS}   & $7$ &  $7$ & $23$ & $23$  & $138$ &  $186$  & $259$  \\
\ctime& \texttt{0.10} & \texttt{0.11} & \texttt{0.32} & \texttt{0.32} & \texttt{1.83} &\texttt{2.38}  & \texttt{3.14} 
  \end{tabular}
  \caption{Performance of POD-Lagrangian suboptimal control for $\hat q = n_t = 18$. BFGS is stopped when $\|\nabla_{\hat u} \hat{\mathcal L} \|_\infty \leq \texttt{tol}_\nabla $, as provided in \cite{Imm}.}
  \label{tab:numbfgsterm}
\end{table}

%% file: plots/1.tikz
% This file was created by matplotlib2tikz v0.5.3.
% The lastest updates can be retrieved from
% 
% https://github.com/nschloe/matplotlib2tikz
% 
% where you can also submit bug reports and leavecomments.
% 
\begin{tikzpicture}

\begin{axis}[
xlabel={x},
ylabel={t},
xmin=0, xmax=1,
ymin=0, ymax=1,
width=\figurewidth,height=\figureheight,axis on top
]
\addplot graphics [includegraphics cmd=\pgfimage,xmin=0, xmax=1, ymin=1, ymax=0] {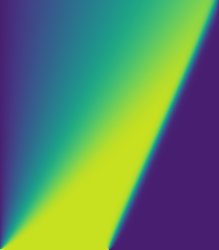};
\end{axis}

\end{tikzpicture}

%% file: plots/2.tikz
% This file was created by matplotlib2tikz v0.5.3.
% The lastest updates can be retrieved from
% 
% https://github.com/nschloe/matplotlib2tikz
% 
% where you can also submit bug reports and leavecomments.
% 
\begin{tikzpicture}

\begin{axis}[
xlabel={x},
ylabel={t},
xmin=0, xmax=1,
ymin=0, ymax=1,
width=\figurewidth,height=\figureheight,axis on top
]
\addplot graphics [includegraphics cmd=\pgfimage,xmin=0, xmax=1, ymin=1, ymax=0] {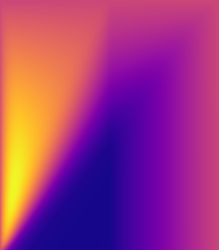};
\end{axis}

\end{tikzpicture}

%% file: plots/3.tikz
% This file was created by matplotlib2tikz v0.5.3.
% The lastest updates can be retrieved from
% 
% https://github.com/nschloe/matplotlib2tikz
% 
% where you can also submit bug reports and leavecomments.
% 
\begin{tikzpicture}

\begin{axis}[
xlabel={x},
ylabel={t},
xmin=0, xmax=1,
ymin=0, ymax=1,
width=\figurewidth,height=\figureheight,axis on top
]
\addplot graphics [includegraphics cmd=\pgfimage,xmin=0, xmax=1, ymin=1, ymax=0] {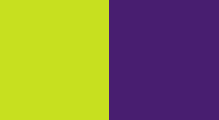};
\end{axis}

\end{tikzpicture}

%% file: plots/4.tikz
% This file was created by matplotlib2tikz v0.5.3.
% The lastest updates can be retrieved from
% 
% https://github.com/nschloe/matplotlib2tikz
% 
% where you can also submit bug reports and leavecomments.
% 
\begin{tikzpicture}

\begin{axis}[
xlabel={x},
ylabel={t},
xmin=0, xmax=1,
ymin=0, ymax=1,
width=\figurewidth,height=\figureheight,axis on top
]
\addplot graphics [includegraphics cmd=\pgfimage,xmin=0, xmax=1, ymin=1, ymax=0] {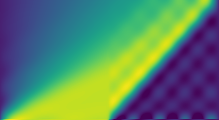};
\end{axis}

\end{tikzpicture}

%% file: plots/5.tikz
% This file was created by matplotlib2tikz v0.5.3.
% The lastest updates can be retrieved from
% 
% https://github.com/nschloe/matplotlib2tikz
% 
% where you can also submit bug reports and leavecomments.
% 
\begin{tikzpicture}

\begin{axis}[
xlabel={x},
ylabel={t},
xmin=0, xmax=1,
ymin=0, ymax=1,
width=\figurewidth,height=\figureheight,axis on top
]
\addplot graphics [includegraphics cmd=\pgfimage,xmin=0, xmax=1, ymin=1, ymax=0] {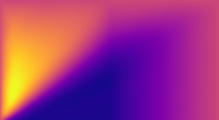};
\end{axis}

\end{tikzpicture}

%% file: plots/6.tikz
% This file was created by matplotlib2tikz v0.5.3.
% The lastest updates can be retrieved from
% 
% https://github.com/nschloe/matplotlib2tikz
% 
% where you can also submit bug reports and leavecomments.
% 
\begin{tikzpicture}

\begin{axis}[
xlabel={x},
ylabel={t},
xmin=0, xmax=1,
ymin=0, ymax=1,
width=\figurewidth,height=\figureheight,axis on top
]
\addplot graphics [includegraphics cmd=\pgfimage,xmin=0, xmax=1, ymin=1, ymax=0] {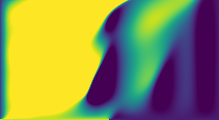};
\end{axis}

\end{tikzpicture}

%% file: plots/7.tikz
% This file was created by matplotlib2tikz v0.5.3.
% The lastest updates can be retrieved from
% 
% https://github.com/nschloe/matplotlib2tikz
% 
% where you can also submit bug reports and leavecomments.
% 
\begin{tikzpicture}

\begin{axis}[
xlabel={x},
ylabel={t},
xmin=0, xmax=1,
ymin=0, ymax=1,
width=\figurewidth,height=\figureheight,axis on top
]
\addplot graphics [includegraphics cmd=\pgfimage,xmin=0, xmax=1, ymin=1, ymax=0] {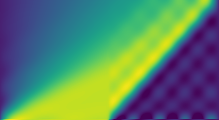};
\end{axis}

\end{tikzpicture}

%% file: plots/8.tikz
% This file was created by matplotlib2tikz v0.5.3.
% The lastest updates can be retrieved from
% 
% https://github.com/nschloe/matplotlib2tikz
% 
% where you can also submit bug reports and leavecomments.
% 
\begin{tikzpicture}

\begin{axis}[
xlabel={x},
ylabel={t},
xmin=0, xmax=1,
ymin=0, ymax=1,
width=\figurewidth,height=\figureheight,axis on top
]
\addplot graphics [includegraphics cmd=\pgfimage,xmin=0, xmax=1, ymin=1, ymax=0] {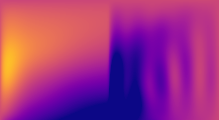};
\end{axis}

\end{tikzpicture}

%% file: plots/9.tikz
% This file was created by matplotlib2tikz v0.5.3.
% The lastest updates can be retrieved from
% 
% https://github.com/nschloe/matplotlib2tikz
% 
% where you can also submit bug reports and leavecomments.
% 
\begin{tikzpicture}

\begin{axis}[
xlabel={x},
ylabel={t},
xmin=0, xmax=1,
ymin=0, ymax=1,
width=\figurewidth,height=\figureheight,axis on top
]
\addplot graphics [includegraphics cmd=\pgfimage,xmin=0, xmax=1, ymin=1, ymax=0] {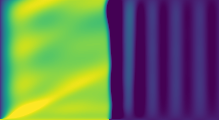};
\end{axis}

\end{tikzpicture}

%% file: plots/10.tikz
% This file was created by matplotlib2tikz v0.5.3.
% The lastest updates can be retrieved from
% 
% https://github.com/nschloe/matplotlib2tikz
% 
% where you can also submit bug reports and leavecomments.
% 
\begin{tikzpicture}

\begin{axis}[
xlabel={x},
ylabel={t},
xmin=0, xmax=1,
ymin=0, ymax=1,
width=\figurewidth,height=\figureheight,axis on top
]
\addplot graphics [includegraphics cmd=\pgfimage,xmin=0, xmax=1, ymin=1, ymax=0] {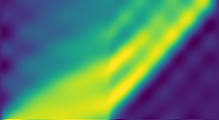};
\end{axis}

\end{tikzpicture}

%% file: plots/11.tikz
% This file was created by matplotlib2tikz v0.5.3.
% The lastest updates can be retrieved from
% 
% https://github.com/nschloe/matplotlib2tikz
% 
% where you can also submit bug reports and leavecomments.
% 
\begin{tikzpicture}

\begin{axis}[
xlabel={x},
ylabel={t},
xmin=0, xmax=1,
ymin=0, ymax=1,
width=\figurewidth,height=\figureheight,axis on top
]
\addplot graphics [includegraphics cmd=\pgfimage,xmin=0, xmax=1, ymin=1, ymax=0] {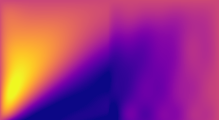};
\end{axis}

\end{tikzpicture}

%% file: plots/12.tikz
% This file was created by matplotlib2tikz v0.5.3.
% The lastest updates can be retrieved from
% 
% https://github.com/nschloe/matplotlib2tikz
% 
% where you can also submit bug reports and leavecomments.
% 
\begin{tikzpicture}

\begin{axis}[
xlabel={x},
ylabel={t},
xmin=0, xmax=1,
ymin=0, ymax=1,
width=\figurewidth,height=\figureheight,axis on top
]
\addplot graphics [includegraphics cmd=\pgfimage,xmin=0, xmax=1, ymin=1, ymax=0] {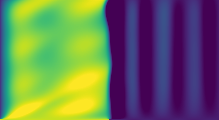};
\end{axis}

\end{tikzpicture}

%% file: chap-conclusion.tex
\section{Conclusion and Outlook}
We have presented a novel approach to low-rank space-time Galerkin approximations that bases on a generalization of classical snapshot-based POD which then can be extended to POD reduction of time discretizations. We have proved optimality of the reduced bases in the relevant function spaces and discussed the numerical implementation.

The space-time Galerkin POD reduction applies well to optimal control problems, as we have illustrated it for the optimal control of a Burgers' equation. Both in terms of computation time for and efficiency of a suboptimal control, the new approach can compete with the combination of classical POD/DEIM for model reduction and BFGS for the optimization. For an optimized distribution of space and time modes, our new approach even dominates the POD/BFGS implementation by achieving better accuracy in less time.

One major potential of the new low-rank space-time Galerkin approach to optimization problems is that it solves the boundary value problem in one shot rather than decoupling forward and backward time like all gradient-based methods do. 

Further possible improvements and issues to future work concerning the proposed space-time POD in application to optimal control problems lye in the freedom of the choice of the measurement functions \cite{BauHS15}. Moreover, the underlying tensor structure is readily extended to include further directions of the state space like parameter dependencies \cite{BauBH15} or inputs. Another issue that needs to be addressed is the treatment of general nonlinearities that can not treated by preassembling like in the presented quadratic case. Then, an inclusion of \emph{empirical interpolation} (EIM) \cite{morBarMNetal04} might be needed to achieve efficiency of the reduction. Moreover, it seems worth investigating whether the principles of space-time POD can be used to construct optimized bases for the interpolation.

\subsection*{Code Availability} \label{sec:codeavailability}
~\\
\begin{minipage}{\linewidth} %minipage to avoid pagebreak inside the framed window
	\begin{framed}
		The source code of the implementations used to compute the presented results can be obtained from:
		\begin{center}
			\href{https://doi.org/10.5281/zenodo.166339}{\texttt{doi:10.5281/zenodo.166339}}
		\end{center}
			
		and is authored by: Manuel Baumann and Jan Heiland \\
		Please contact Manuel Baumann and Jan Heiland for licensing information
	\end{framed}
\end{minipage}

\subsection*{Acknowledgements} 
We thank Joost van Zwieten, co-developer of \textit{Nutils}\footnote{Open source finite element toolbox for Python: \href{http://nutils.org}{http://nutils.org}}, for providing benchmarks and valuable insight into space-time discretizations of Burgers' equation.